\documentclass[12pt,thmsa]{article}
\usepackage{amssymb}
\usepackage{amsfonts}
\usepackage{sw20bams}

\input tcilatex
\begin{document}

\author{Steven R. Finch}
\title{Oblique Circular Cones and Cylinders}
\date{December 31, 2012}
\maketitle

\begin{abstract}
Surface area and mean width of a cylinder (the convex hull of two parallel
disks) in $\mathbb{R}^{3}$ are computed. It is more difficult to obtain
analogous results for a cone (the convex hull of a disk $D$ and a point $p$%
). \ Oblique formulas for mean width, as well as those for mean curvature,
are new. Let $\ell $ denote the unique diameter of $D$ whose endpoints are
equidistant from $p$. We conclude with a question involving the plane that
bisects the cone and contains $\{p,\ell \}$, as $p$ varies. What is the
minimum ratio of the smaller measure to the larger? \ 
\end{abstract}

\footnotetext{%
Copyright \copyright\ 2012 by Steven R. Finch. All rights reserved.}Certain
mathematical formulas need better marketing. An example is the formula for
the surface area of an oblique circular cone. \ Although an exact expression
has been known since at least 1825, it has been rediscovered several times
by independent researchers. \ The same could be said for the surface area of
an oblique circular cylinder, but the redundancy of effort here is less
costly (the formula follows swiftly from a well-known expression for the
circumference of an ellipse).

This paper continues our work \cite{Fi1} on convex hulls of two disks in $%
\mathbb{R}^{3}$. \ Three numerical characteristics of a convex hull --
volume, surface area and mean width -- are central. These quantities, along
with the Euler characteristic, form a basis of the space of all additive
continuous measures that are invariant under rigid motions. \
\textquotedblleft The mean width is a new measure on three-dimensional
solids that enjoys equal rights with volume and surface
area\textquotedblright\ \cite{Rt}. As far as we know, our expressions for
mean width in the oblique case have not appeared before.

If we split an oblique cylinder along its axis, the two resultant
half-cylinders have equal $VL$, $AR$ and $MW$. \ The same symmetry does not
hold for an oblique cone. \ While $VL_{1}=VL_{2}$ for the half-cones, the
associations between $AR_{1}\leq AR_{2}$ and $MW_{1}\leq MW_{2}$ are more
complicated. \ An optimization problem involving the ratios $AR_{1}/AR_{2}$
and $MW_{1}/MW_{2}$ will occupy us at the end. \ We agree to split a cone so
that its cross-section is an isosceles triangle (uniquely); a removal of
this requirement is left open for someone else to address. \ There is an
issue of whether we employ total surface area (including bases) or lateral
surface area (curved portions only). \ The parameters minimizing $%
AR_{1}/AR_{2}$ and $MW_{1}/MW_{2}$ are nontrivial, as will be seen.

\section{Oblique Cylinder}

The solid cylinder $\Omega $ is the convex hull of the following two
parallel disks in $\mathbb{R}^{3}$:%
\[
\begin{array}{ccccc}
\left\{ (x,y,z):x^{2}+y^{2}\leq 1\text{ \&\ }z=0\right\} &  & \text{and} & 
& \left\{ (x,y,z):x^{2}+(y-a)^{2}\leq 1\text{ \&\ }z=b\right\}%
\end{array}%
\]%
where WLOG $a\geq 0$, $b>0$. It is easy to find the volume of $\Omega $
(integrating the area of horizontal slices):%
\[
VL=\dint\limits_{0}^{b}\pi \,dz=\pi b. 
\]%
The curved portion of the boundary $\partial \Omega $ in the half-space $%
x\geq 0$ can be represented parametrically:%
\[
\begin{array}{ccccccc}
x=\sqrt{1-v^{2}}, &  & y=au+v, &  & z=bu, &  & 0\leq u\leq 1\text{ \&\ }%
-1\leq v\leq 1.%
\end{array}%
\]%
Let $x_{u}$, $x_{v}$, $y_{u}$, $y_{v}$, $z_{u}$, $z_{v}$ denote partial
derivatives of $x$, $y$, $z$. \ Defining%
\[
\begin{array}{ccc}
E=\left( x_{u},y_{u},z_{u}\right) \cdot \left( x_{u},y_{u},z_{u}\right) , & 
& G=\left( x_{v},y_{v},z_{v}\right) \cdot \left( x_{v},y_{v},z_{v}\right) ,%
\end{array}%
\]%
\[
F=\left( x_{u},y_{u},z_{u}\right) \cdot \left( x_{v},y_{v},z_{v}\right) 
\]%
we have%
\begin{eqnarray*}
AR &=&\dint\limits_{\partial \Omega }dS=2\pi
+2\dint\limits_{0}^{1}\dint\limits_{-1}^{1}\sqrt{EG-F^{2}}\,dv\,du \\
&=&2\pi +2\dint\limits_{0}^{1}\dint\limits_{-1}^{1}\sqrt{\frac{%
a^{2}v^{2}+b^{2}}{1-v^{2}}}\,dv\,du \\
&=&2\left( \pi +2\sqrt{a^{2}+b^{2}}E\left( \dfrac{a^{2}}{a^{2}+b^{2}}\right)
\right)
\end{eqnarray*}%
where 
\[
E(\mu )=\dint\limits_{0}^{\pi /2}\sqrt{1-\mu \sin (\theta )^{2}}\,d\theta
=\dint\limits_{0}^{1}\sqrt{\dfrac{1-\mu \,t^{2}}{1-t^{2}}}\,dt 
\]%
is the complete elliptic integral of the second kind.

To compute $MW$ (via the \textquotedblleft indirect
approach\textquotedblright\ in \cite{Fi1}), first find the mean curvature $H$
and integrate. Defining%
\[
\mathcal{N}=-\frac{\left( x_{u},y_{u},z_{u}\right) \times \left(
x_{v},y_{v},z_{v}\right) }{\left\vert \left( x_{u},y_{u},z_{u}\right) \times
\left( x_{v},y_{v},z_{v}\right) \right\vert }=\left( b\sqrt{\frac{1-v^{2}}{%
a^{2}v^{2}+b^{2}}},\frac{bv}{\sqrt{a^{2}v^{2}+b^{2}}},\frac{-av}{\sqrt{%
a^{2}v^{2}+b^{2}}}\right) , 
\]%
\[
\begin{array}{ccc}
L=\left( x_{u},y_{u},z_{u}\right) \cdot \mathcal{N}_{u}, &  & N=\left(
x_{v},y_{v},z_{v}\right) \cdot \mathcal{N}_{v},%
\end{array}%
\]%
\[
M=\dfrac{1}{2}\left( \left( x_{u},y_{u},z_{u}\right) \cdot \mathcal{N}%
_{v}+\left( x_{v},y_{v},z_{v}\right) \cdot \mathcal{N}_{u}\right) 
\]%
we have%
\begin{eqnarray*}
\dint\limits_{\partial \Omega }H\,dS &=&\dint\limits_{\partial \Omega }\frac{%
EN-2FM+GL}{2\left( EG-F^{2}\right) }\,dS \\
&=&2\dint\limits_{0}^{1}\dint\limits_{-1}^{1}\frac{EN-2FM+GL}{2\left(
EG-F^{2}\right) }\sqrt{EG-F^{2}}\,dv\,du \\
&=&2\dint\limits_{0}^{1}\dint\limits_{-1}^{1}\frac{\left( a^{2}+b^{2}\right)
b}{2\left( a^{2}v^{2}+b^{2}\right) ^{3/2}}\sqrt{\frac{a^{2}v^{2}+b^{2}}{%
1-v^{2}}}\,dv\,du \\
&=&\dint\limits_{-1}^{1}\frac{\left( a^{2}+b^{2}\right) b}{a^{2}v^{2}+b^{2}}%
\,\frac{1}{\sqrt{1-v^{2}}}\,dv=\pi \sqrt{a^{2}+b^{2}}.
\end{eqnarray*}%
Second, find the exterior dihedral angle along each edge and integrate. \
For example, the semicircular edge $x^{2}+y^{2}=1$, $z=0$ \&\ $x\geq 0$
corresponds to $u=0$ \&\ $-1\leq v\leq 1$. \ Call this $\varepsilon $. \ The
exterior dihedral angle is%
\[
\alpha =\arccos \left( \mathcal{N}\cdot (0,0,-1)\right) =\arccos \left( 
\frac{av}{\sqrt{a^{2}v^{2}+b^{2}}}\right) 
\]%
and arclength $s$ satisfies%
\[
ds=\sqrt{x_{v}^{2}+y_{v}^{2}}\,dv=\frac{1}{\sqrt{1-v^{2}}}\,dv. 
\]%
It follows that 
\[
\dint\limits_{\varepsilon }\alpha \,ds=\dint\limits_{-1}^{1}\frac{1}{\sqrt{%
1-v^{2}}}\arccos \left( \frac{av}{\sqrt{a^{2}v^{2}+b^{2}}}\right) dv=\frac{1%
}{2}\pi ^{2}. 
\]%
Finally, the surface $\partial \Omega $ is piecewise continuously
differentiable and has $n=4$ smooth edges $\varepsilon _{j}$ with
(non-constant)\ dihedral angles $\alpha _{j}$, $1\leq j\leq n$. From the
general formula 
\[
MW=\frac{1}{2\pi }\dint\limits_{\partial \Omega }H\,dS+\frac{1}{4\pi }%
\sum_{j=1}^{n}\dint\limits_{\varepsilon _{j}}\alpha _{j}\,ds, 
\]%
we deduce that%
\[
MW=\frac{1}{2\pi }\left( \pi \sqrt{a^{2}+b^{2}}\right) +\frac{1}{4\pi }%
\left( \frac{4}{2}\pi ^{2}\right) =\frac{1}{2}\left( \sqrt{a^{2}+b^{2}}+\pi
\right) . 
\]%
In the special case $a=0$ (a right circular cylinder), $MW=\left( b+\pi
\right) /2$, which is consistent with \cite{SKM, Sa}.

\subsection{Comment about Ellipses}

It is well known \cite{SL, Ry, Ptr} that the lateral surface area of a
cylinder is the product of an element length, $\sqrt{a^{2}+b^{2}}$, and the
circumference of a perpendicular section. The latter is an ellipse with
semi-major axis $1$ and semi-minor axis $b/\sqrt{a^{2}+b^{2}}$. \ This is
the shadow cast by the unit $xy$-circle into any plane (in $xyz$-space)
normal to the vector $(0,a,b)$. \ In particular, the length of the
projection of vector $(0,1,0)$ onto vector $(0,b,-a)$ is 
\[
\frac{(0,1,0)\cdot (0,b,-a)}{\sqrt{a^{2}+b^{2}}}=\frac{b}{\sqrt{a^{2}+b^{2}}}
\]%
which is the desired semi-minor axis. The eccentricity squared is hence%
\[
e^{2}=1-\frac{b^{2}}{a^{2}+b^{2}}=\frac{a^{2}}{a^{2}+b^{2}} 
\]%
and the circumference of the ellipse is $4E(e^{2})$. \ 

Older references discussing the surface area of a cylinder include \cite{Rs,
Ly, Bk}. \ A new result \cite{Adlj} provides an AGM-like iteration that
permits rapid calculation of $E(e^{2})$.

\section{Oblique Cone}

The solid cone $\Omega $ is the convex hull of a disk and a point in $%
\mathbb{R}^{3}$:%
\[
\begin{array}{ccccc}
\left\{ (x,y,z):x^{2}+y^{2}\leq 1\text{ \&\ }z=0\right\} &  & \text{and} & 
& (0,a,b)%
\end{array}%
\]%
where WLOG $a\geq 0$, $b>0$. It is easy, as before, to find the volume of $%
\Omega $:%
\[
VL=\dint\limits_{0}^{b}\pi \left( 1-\frac{z}{b}\right) ^{2}dz=\left. -\frac{1%
}{3}\pi \left( 1-\frac{z}{b}\right) ^{3}b\right\vert _{0}^{b}=\frac{1}{3}\pi
b. 
\]%
The curved portion of the boundary $\partial \Omega $ in the half-space $%
x\geq 0$ can be represented parametrically:%
\[
\begin{array}{ccccccc}
x=(1-u)\sqrt{1-v^{2}}, &  & y=au+(1-u)v, &  & z=bu, &  & 0\leq u\leq 1\text{
\&\ }-1\leq v\leq 1.%
\end{array}%
\]%
After computing $E$, $F$, $G$, we have%
\begin{eqnarray*}
AR &=&\dint\limits_{\partial \Omega }dS=\pi
+2\dint\limits_{0}^{1}\dint\limits_{-1}^{1}\sqrt{EG-F^{2}}\,dv\,du \\
&=&\pi +2\dint\limits_{0}^{1}\dint\limits_{-1}^{1}(1-u)\sqrt{\frac{%
(1-av)^{2}+b^{2}}{1-v^{2}}}\,dv\,du \\
&=&\pi +2\sqrt{s_{0}s_{1}}\left[ E(c_{0})-K(c_{0})+(1-c_{1})\Pi (c_{1},c_{0})%
\right]
\end{eqnarray*}%
where%
\[
\begin{array}{ccc}
s_{0}=\sqrt{(1-a)^{2}+b^{2}}, &  & s_{1}=\sqrt{(1+a)^{2}+b^{2}},%
\end{array}%
\]%
\[
\begin{array}{ccc}
c_{0}=\dfrac{1}{2}\left( 1-\dfrac{1-a^{2}+b^{2}}{s_{0}s_{1}}\right) , &  & 
c_{1}=\dfrac{1}{2}\left( 1-\dfrac{1+a^{2}+b^{2}}{s_{0}s_{1}}\right) .%
\end{array}%
\]%
In the preceding,%
\[
K(\mu )=\dint\limits_{0}^{\pi /2}\dfrac{1}{\sqrt{1-\mu \sin (\theta )^{2}}}%
\,d\theta =\dint\limits_{0}^{1}\dfrac{1}{\sqrt{(1-t^{2})(1-\mu \,t^{2})}}%
\,dt 
\]%
is the complete elliptic integral of the first kind and 
\[
\Pi (\nu ,\mu )=\dint\limits_{0}^{\pi /2}\dfrac{1}{\left( 1-\nu \sin (\theta
)^{2}\right) \sqrt{1-\mu \sin (\theta )^{2}}}\,d\theta =\dint\limits_{0}^{1}%
\dfrac{1}{(1-\nu \,t^{2})\sqrt{(1-t^{2})(1-\mu \,t^{2})}}\,dt 
\]%
is the complete elliptic integral of the third kind. \ Interestingly, a
rapid AGM-like iteration to calculate $\Pi (c_{1},c_{0})$ would seem still
to be awaiting discovery. \ This would \textquotedblleft finish the
quest\textquotedblright\ because such iterations for $E(c_{0})$ and $%
K(c_{0}) $ are already in our possession \cite{Adlj}. \ 

As before,%
\begin{eqnarray*}
\mathcal{N} &=&-\frac{\left( x_{u},y_{u},z_{u}\right) \times \left(
x_{v},y_{v},z_{v}\right) }{\left\vert \left( x_{u},y_{u},z_{u}\right) \times
\left( x_{v},y_{v},z_{v}\right) \right\vert } \\
&=&\left( b\sqrt{\frac{1-v^{2}}{(1-av)^{2}+b^{2}}},\frac{bv}{\sqrt{%
(1-av)^{2}+b^{2}}},\frac{1-av}{\sqrt{(1-av)^{2}+b^{2}}}\right)
\end{eqnarray*}%
and%
\begin{eqnarray*}
\dint\limits_{\partial \Omega }H\,dS
&=&2\dint\limits_{0}^{1}\dint\limits_{-1}^{-1}\frac{EN-2FM+GL}{2\left(
EG-F^{2}\right) }\sqrt{EG-F^{2}}\,dv\,du \\
&=&2\dint\limits_{0}^{1}\dint\limits_{-1}^{-1}\frac{(1+a^{2}+b^{2}-2av)b}{%
2(1-u)\left( (1-av)^{2}+b^{2}\right) ^{3/2}} \\
&&\cdot \,(1-u)\sqrt{\frac{(1-av)^{2}+b^{2}}{1-v^{2}}}\,dv\,du \\
&=&\dint\limits_{-1}^{1}\frac{(1+a^{2}+b^{2}-2av)b}{(1-av)^{2}+b^{2}}\frac{1%
}{\sqrt{1-v^{2}}}\,dv \\
&=&\frac{1}{2}\left( \sqrt{a^{2}+(-i+b)^{2}}+\sqrt{a^{2}+(i+b)^{2}}\right)
\pi
\end{eqnarray*}%
where $i$ is the imaginary unit and we take the branch cut along the
negative real axis (for both square root and logarithm functions). \ 

Let $\varepsilon $ denote the arc of the semicircle $x^{2}+y^{2}=1$, $z=0$
\&\ $x\geq 0$. \ The exterior dihedral angle is%
\[
\alpha =\arccos \left( \mathcal{N}\cdot (0,0,-1)\right) =\arccos \left( 
\frac{-1+av}{\sqrt{(1-av)^{2}+b^{2}}}\right) 
\]%
and arclength $s$ satisfies%
\[
ds=\sqrt{x_{v}^{2}+y_{v}^{2}}\,dv=\frac{1}{\sqrt{1-v^{2}}}\,dv. 
\]%
It follows that 
\begin{eqnarray*}
\dint\limits_{\varepsilon }\alpha \,ds &=&\dint\limits_{-1}^{1}\frac{1}{%
\sqrt{1-v^{2}}}\arccos \left( \frac{-1+av}{\sqrt{(1-av)^{2}+b^{2}}}\right) dv
\\
&=&\frac{1}{2}\left[ \pi +i\,\ln \left( -i+b+\sqrt{a^{2}+(-i+b)^{2}}\right)
-i\,\ln \left( i+b+\sqrt{a^{2}+(i+b)^{2}}\right) \right] \pi
\end{eqnarray*}%
and therefore \ 
\begin{eqnarray*}
MW &=&\frac{1}{2\pi }\dint\limits_{\partial \Omega }H\,dS+\frac{2}{4\pi }%
\dint\limits_{\varepsilon }\alpha \,ds \\
&=&\frac{1}{4}\left( \sqrt{a^{2}+(-i+b)^{2}}+\sqrt{a^{2}+(i+b)^{2}}\right) \\
&&+\,\frac{1}{4}\left[ \pi +i\,\ln \left( -i+b+\sqrt{a^{2}+(-i+b)^{2}}%
\right) -i\,\ln \left( i+b+\sqrt{a^{2}+(i+b)^{2}}\right) \right] .
\end{eqnarray*}%
In the special case $a=0$ (a right circular cone),%
\begin{eqnarray*}
MW &=&\frac{1}{4}\left[ (-i+b)+(i+b)\right] +\frac{1}{4}\left[ \pi +i\,\ln
\left( -2i+2b\right) -i\,\ln \left( 2i+2b\right) \right] \\
&=&\frac{b}{2}+\frac{\pi }{4}+\frac{i}{4}\left[ \ln (2)+\ln \left(
i-b\right) -i\pi \right] -\frac{i}{4}\left[ \ln (2)+\ln \left( i+b\right) %
\right] \\
&=&\frac{b}{2}+\frac{\pi }{2}-\frac{1}{2}\left[ \frac{i}{2}\ln \left(
i+b\right) -\frac{i}{2}\ln \left( i-b\right) \right] =\frac{b+\pi }{2}-\frac{%
1}{2}\arctan (b)
\end{eqnarray*}%
which is consistent with \cite{Sa}.

\subsection{Comment about Ruled Surfaces}

The curved portion of a cone is a ruled surface, which means that through
every point there exists a linear element that lies on the surface. \ A
mistaken argument in \cite{Mz} gave that the lateral surface area of a right
circular cone is $2\pi \sqrt{1+b^{2}}$, which seems natural since the
element length is $\sqrt{1+b^{2}}$ and the base circumference is $2\pi $. \
This is an error:\ the correct value \cite{SL} should be $\pi \sqrt{1+b^{2}}$
(seen by setting $a=0$, implying that $c_{0}=c_{1}=0$, in our formula for $%
AR $). Alternatively, we might use Pappus's centroid theorem or a
solid-of-revolution approach: letting $y=1-z/b$, 
\begin{eqnarray*}
AR &=&\dint\limits_{0}^{b}2\pi y\sqrt{1+\left( \frac{dy}{dz}\right) ^{2}}%
dz=2\pi \sqrt{1+\frac{1}{b^{2}}}\dint\limits_{0}^{b}\left( 1-\frac{z}{b}%
\right) dz \\
&=&2\pi \sqrt{\frac{1}{b^{2}}+1}\cdot \left. \left( -\frac{b}{2}\right)
\left( 1-\frac{z}{b}\right) ^{2}\right\vert _{0}^{b}=\pi \sqrt{1+b^{2}}.
\end{eqnarray*}%
The mistake was generalized to oblique circular cones in \cite{Mz}, under
the hope that integrating the (varying) element length would provide the
surface area. This too is an error: the correct procedure is found in \cite%
{Ry}. Also, the double angle formula for cosine was misapplied, hence the
last coefficient \textquotedblleft 2\textquotedblright\ in formula (6) of 
\cite{Mz} should be \textquotedblleft 4\textquotedblright\ (and likewise in
subsequent formulas).

Legendre \cite{Lg} was the first person to explicitly give $AR$ in terms of
elliptic integrals. \ Follow-up work appeared in \cite{Hu, Tr, Sc, Fu, Ws,
Lb, Fp, Pw, Cn}. Older references discussing the surface area of a cone
include \cite{Ly, Ms, Hr, Itd, Eu, Kr}. The oblique formula for $AR$ was
publicized early in the \textit{Edinburgh Encyclopedia }\cite{Wc}, \textit{%
Encyclopedia Britannica} \cite{Wi} and elsewhere \cite{Hy, Mr, Gh, Pra};
thus the marketing fiasco has occurred only relatively recently. \
Legendre's achievement deserves not to be forgotten!

\section{Optimization Preliminaries}

Considering all cones $\Omega $ of equal volume, which one has the least
surface area? An elementary solution (without use of elliptic integrals) is
obtained by differentiating 
\[
AR-\pi =\dint\limits_{-1}^{1}\sqrt{\frac{(1-av)^{2}+b^{2}}{1-v^{2}}}\,dv 
\]%
twice with respect to $a$:%
\[
\frac{\partial }{\partial a}\dint\limits_{-1}^{1}\sqrt{\frac{(1-av)^{2}+b^{2}%
}{1-v^{2}}}\,dv=\dint\limits_{-1}^{1}\frac{-v(1-av)}{\sqrt{(1-av)^{2}+b^{2}}}%
\frac{1}{\sqrt{1-v^{2}}}\,dv, 
\]%
\[
\frac{\partial ^{2}}{\partial a^{2}}\dint\limits_{-1}^{1}\sqrt{\frac{%
(1-av)^{2}+b^{2}}{1-v^{2}}}\,dv=\dint\limits_{-1}^{1}\frac{b^{2}v^{2}}{\left[
(1-av)^{2}+b^{2}\right] ^{3/2}}\frac{1}{\sqrt{1-v^{2}}}\,dv. 
\]%
Setting $a=0$ in the former expression, an odd function appears inside the
integral, hence the first derivative vanishes. \ The second derivative is
always positive, ensuring that the right circular cone is minimal. \ Our
approach follows \cite{DF1, DF2}, but caution should be exercised: beginning
with formula (2.1) of \cite{DF1},\ $a^{2}+h^{2}+r^{2}-2ar\cos (\theta )$
should be everywhere replaced by $a^{2}\cos (\theta )^{2}+h+r^{2}-2ar\cos
(\theta )$. The same error occurs in \cite{DF2}. Again, see \cite{Ry} for
supporting details.

Considering all cones $\Omega $ of equal volume, which one has the least
mean width? Starting with the algebraic expression for $MW$, it is easily
verified that the right circular cone is minimal here too.

Before turning to more complicated optimization problems, let us define $\xi
(x)$ for any $x>0$ to be the largest real zero $y$ of $xy-\cos (y)$. \ The
most famous value of this function is%
\[
\xi (1)=0.7390851332151606416553120..., 
\]%
popularized in \cite{Kap}. \ Define also $\eta (x)$ for any $0<x<1$ to be
the largest real zero $y$ of $xy-\sin (y)$. \ We will require two values:%
\[
\xi \left( \frac{2}{\pi }\right) =0.9340137863539518545607006..., 
\]%
\[
\eta \left( \frac{1}{\pi }\right) =2.3137341320786811322489898.... 
\]

In preparing the following solutions, we are being somewhat hasty. Rigorous
proofs of minimality are not found; high-precision numerical confirmations
often constitute the basis of our intuition.

\section{Ratio of Surface Areas}

Exact expressions for the two integrals:%
\[
\begin{array}{ccccc}
\dint\limits_{-1}^{0}\sqrt{\dfrac{(1-av)^{2}+b^{2}}{1-v^{2}}}\,dv &  & \text{%
and} &  & \dint\limits_{0}^{1}\sqrt{\dfrac{(1-av)^{2}+b^{2}}{1-v^{2}}}\,dv%
\end{array}%
\]
seem to be infeasible, despite the fact that we know already their \textit{%
sum} and the availability of various transformations \cite{BF}. \ This will
not be a deterrent, however. because only the case $a>1$ \& $b\rightarrow
0^{+}$ is conjectured to be relevant.

\subsection{Lateral AR}

Since $a>1$, we have%
\begin{eqnarray*}
\lim_{b\rightarrow 0^{+}}\frac{\dint\limits_{0}^{1}\sqrt{\dfrac{%
(1-av)^{2}+b^{2}}{1-v^{2}}}\,dv}{\dint\limits_{-1}^{0}\sqrt{\dfrac{%
(1-av)^{2}+b^{2}}{1-v^{2}}}\,dv} &=&\frac{\dint\limits_{0}^{1/a}\dfrac{1-av}{%
\sqrt{1-v^{2}}}\,dv-\dint\limits_{1/a}^{1}\dfrac{1-av}{\sqrt{1-v^{2}}}\,dv}{%
\dint\limits_{-1}^{0}\dfrac{1-av}{\sqrt{1-v^{2}}}\,dv} \\
&=&\frac{\pi -2a+4\sqrt{a^{2}-1}-4\func{arcsec}(a)}{\pi +2a}.
\end{eqnarray*}%
After differentiating, the equation to be solved is%
\[
\frac{\pi }{2a}\sqrt{a^{2}-1}=\func{arccsc}(a) 
\]%
hence%
\[
a=1.2437608987462040336147443.... 
\]%
Further,%
\[
\sin \left( \frac{\pi }{2a}\sqrt{a^{2}-1}\right) =\frac{1}{a} 
\]%
which (after defining $x$ to be the argument of $\sin $) implies%
\[
\frac{\pi }{\sqrt{\pi ^{2}-4x^{2}}}\sin (x)=1 
\]%
which implies%
\[
\pi ^{2}\sin (x)^{2}=\pi ^{2}-4x^{2} 
\]%
which implies%
\[
4x^{2}=\pi ^{2}\cos (x)^{2} 
\]%
which implies%
\[
2x=\pi \cos (x) 
\]%
hence $x=\xi \left( 2/\pi \right) $. \ Thus the infimum of ratios is $%
0.1892....$

\subsection{Total AR}

Since $a>1$, we have%
\begin{eqnarray*}
\lim_{b\rightarrow 0^{+}}\frac{\dfrac{\pi }{2}+\dint\limits_{0}^{1}\sqrt{%
\dfrac{(1-av)^{2}+b^{2}}{1-v^{2}}}\,dv}{\dfrac{\pi }{2}+\dint\limits_{-1}^{0}%
\sqrt{\dfrac{(1-av)^{2}+b^{2}}{1-v^{2}}}\,dv} &=&\frac{\dfrac{\pi }{2}%
+\dint\limits_{0}^{1/a}\dfrac{1-av}{\sqrt{1-v^{2}}}\,dv-\dint%
\limits_{1/a}^{1}\dfrac{1-av}{\sqrt{1-v^{2}}}\,dv}{\dfrac{\pi }{2}%
+\dint\limits_{-1}^{0}\dfrac{1-av}{\sqrt{1-v^{2}}}\,dv} \\
&=&\frac{-a+2\sqrt{a^{2}-1}+2\func{arccsc}(a)}{\pi +a}.
\end{eqnarray*}%
After differentiating, the equation to be solved is%
\[
\frac{\pi }{a}\left( a-\sqrt{a^{2}-1}\right) =\func{arcsec}(a) 
\]%
hence%
\[
a=1.4782960807222794430758369.... 
\]%
Further,%
\[
\cos \left( \frac{\pi }{a}\left( a-\sqrt{a^{2}-1}\right) \right) =\frac{1}{a}
\]%
which (after defining $x$ to be the argument of $\cos $) implies%
\[
\frac{\pi }{\sqrt{2\pi x-x^{2}}}\cos (x)=1 
\]%
which implies%
\[
\pi ^{2}\cos (x)^{2}=2\pi x-x^{2} 
\]%
which implies%
\[
(\pi -x)^{2}=\pi ^{2}-2\pi x+x^{2}=\pi ^{2}-\pi ^{2}\cos (x)^{2}=\pi
^{2}\sin (x)^{2} 
\]%
which implies%
\[
\pi -x=\pi \sin (x) 
\]%
which (after defining $w=\pi -x$) implies%
\[
w=\pi \sin (\pi -w)=\pi \sin (w) 
\]%
hence $w=\eta \left( 1/\pi \right) $. Thus the infimum of ratios is $%
0.4729....$

\section{Ratio of Mean Widths}

Unlike the preceding, exact expressions for the two integrals:%
\[
J_{0}=\dint\limits_{-1}^{0}\dfrac{(1+a^{2}+b^{2}-2av)b}{(1-av)^{2}+b^{2}}%
\dfrac{1}{\sqrt{1-v^{2}}}\,dv, 
\]%
\[
J_{1}=\dint\limits_{0}^{1}\dfrac{(1+a^{2}+b^{2}-2av)b}{(1-av)^{2}+b^{2}}%
\dfrac{1}{\sqrt{1-v^{2}}}\,dv 
\]%
are needed. \ Only the case $a>1$ \& $b\rightarrow 0^{+}$ is conjectured to
be relevant. \ Integral $J_{1}$ is problematic due to the singularity at $%
v=1/a$ as $b\rightarrow 0^{+}$, whereas integral $J_{0}\rightarrow 0$
clearly.

\subsection{Integrated Mean Curvature}

It is not difficult to show that%
\begin{eqnarray*}
J_{0} &=&\frac{i}{2}\left\{ \sqrt{a^{2}+(-i+b)^{2}}\left[ \frac{i\pi }{2}%
-\ln (i-b)+\ln \left( a+\sqrt{a^{2}+(-i+b)^{2}}\right) \right] \right. \\
&&+\,\left. \sqrt{a^{2}+(i+b)^{2}}\left[ \frac{i\pi }{2}+\ln (-i-b)-\ln
\left( a+\sqrt{a^{2}+(i+b)^{2}}\right) \right] \right\}
\end{eqnarray*}%
From the known expression for $J_{0}+J_{1}$, we deduce that%
\begin{eqnarray*}
J_{1} &=&\frac{i}{2}\left\{ \sqrt{a^{2}+(-i+b)^{2}}\left[ -\frac{3i\pi }{2}%
+\ln (i-b)-\ln \left( a+\sqrt{a^{2}+(-i+b)^{2}}\right) \right] \right. \\
&&+\,\left. \sqrt{a^{2}+(i+b)^{2}}\left[ -\frac{3i\pi }{2}-\ln (-i-b)+\ln
\left( a+\sqrt{a^{2}+(i+b)^{2}}\right) \right] \right\}
\end{eqnarray*}%
and therefore $J_{1}\rightarrow \pi \sqrt{a^{2}-1}$.

\subsection{Exterior Dihedral Angles}

The isosceles triangle contains vertices $(0,a,b)$, $(-1,0,0)$, $(1,0,0)$,
hence is contained in the plane $-by+az=0$. \ The exterior normal vector to
this face depends on the choice of half-cone. The base of the triangle is of
length $2$; the sum of the lengths of the two triangular legs is $2\sqrt{%
a^{2}+b^{2}+1}$.

For the lower (smaller) half-cone, the exterior normal vector is $(0,-b,a)$.
At the base edge, the dihedral angle is 
\[
\arccos \left( (0,0,-1)\cdot \frac{(0,-b,a)}{\sqrt{a^{2}+b^{2}}}\right) =\pi
-\arccos \left( \frac{a}{\sqrt{a^{2}+b^{2}}}\right) 
\]%
and at each leg edge (for which $v=0$), the dihedral angle is

\[
\arccos \left( \mathcal{N}\cdot \frac{(0,-b,a)}{\sqrt{a^{2}+b^{2}}}\right)
=\arccos \left( \frac{a}{\sqrt{1+b^{2}}\sqrt{a^{2}+b^{2}}}\right) . 
\]

For the upper (larger) half-cone, the exterior normal vector is $(0,b,-a)$.
At the base edge, the dihedral angle is 
\[
\arccos \left( (0,0,-1)\cdot \frac{(0,b,-a)}{\sqrt{a^{2}+b^{2}}}\right)
=\arccos \left( \frac{a}{\sqrt{a^{2}+b^{2}}}\right) 
\]%
and at each leg edge, the dihedral angle is

\[
\arccos \left( \mathcal{N}\cdot \frac{(0,b,-a)}{\sqrt{a^{2}+b^{2}}}\right)
=\pi -\arccos \left( \frac{a}{\sqrt{1+b^{2}}\sqrt{a^{2}+b^{2}}}\right) . 
\]%
\qquad

Note also that 
\[
L_{0}=2\dint\limits_{-1}^{0}\dfrac{1}{\sqrt{1-v^{2}}}\arccos \left( \dfrac{%
-1+av}{\sqrt{(1-av)^{2}+b^{2}}}\right) dv\rightarrow 2\dint\limits_{-1}^{0}%
\dfrac{\pi }{\sqrt{1-v^{2}}}\,dv=\pi ^{2}, 
\]%
\[
L_{1}=2\dint\limits_{0}^{1}\dfrac{1}{\sqrt{1-v^{2}}}\arccos \left( \dfrac{%
-1+av}{\sqrt{(1-av)^{2}+b^{2}}}\right) dv\rightarrow 2\dint\limits_{0}^{1/a}%
\dfrac{\pi }{\sqrt{1-v^{2}}}\,dv=2\pi \func{arccsc}(a) 
\]
as $b\rightarrow 0^{+}$, for upper/lower semicircular edges.

\subsection{Culmination}

As $b\rightarrow 0^{+}$, the ratio 
\[
\frac{\dfrac{J_{1}}{2\pi }+\dfrac{1}{4\pi }\left[ L_{1}+2\left( \pi -\arccos
\left( \dfrac{a}{\sqrt{a^{2}+b^{2}}}\right) \right) +2\sqrt{a^{2}+b^{2}+1}%
\arccos \left( \dfrac{a}{\sqrt{1+b^{2}}\sqrt{a^{2}+b^{2}}}\right) \right] }{%
\dfrac{J_{0}}{2\pi }+\dfrac{1}{4\pi }\left[ L_{0}+2\arccos \left( \dfrac{a}{%
\sqrt{a^{2}+b^{2}}}\right) +2\sqrt{a^{2}+b^{2}+1}\left( \pi -\arccos \left( 
\dfrac{a}{\sqrt{1+b^{2}}\sqrt{a^{2}+b^{2}}}\right) \right) \right] } 
\]%
approaches%
\[
\frac{\dfrac{1}{2}\sqrt{a^{2}-1}+\dfrac{1}{2}\func{arccsc}(a)+\dfrac{1}{2}}{%
\dfrac{\pi }{4}+\dfrac{1}{2}\sqrt{a^{2}+1}}. 
\]%
After differentiating, the equation to be solved is%
\[
\frac{\sqrt{a^{2}-1}\left( 2+\pi \sqrt{a^{2}+1}\right) }{2a^{2}}=1+\func{%
arccsc}(a) 
\]%
hence%
\[
a=1.3638337555895594010839152.... 
\]%
Thus the infimum of ratios is $0.8431....$

\section{Acknowledgements}

Wouter Meeussen's package ConvexHull3D.m was helpful to me in preparing this
paper \cite{Msn}. \ He kindly extended the software functionality at my
request. \ I\ would be grateful for assistance in expanding my bibliography:
surely I have missed more than a few documents on surface area (for oblique
circular cylinders and cones)!

\section{Addendum}

There is a curious asymmetry in our discussions of $AR$ in Section 4 and of $%
MW$\ in Section 5. For $AR$, we imagine a fixed amount of sheet metal,
allocated between the two half-cones (thus $AR_{1}+AR_{2}$ is equal to the
original $AR$). In contrast, for $MW$, we treat the half-cones as
independent convex solids -- think of widths as drawn from wholly separate
objects -- although each object contains the isosceles triangle as a face
(thus $MW_{1}+MW_{2}$ is larger than the original $MW$).

To make things compatible, let us revisit Section 4.2, but now include the
overlap. \ The triangular area is $\frac{1}{2}\cdot 2\cdot \sqrt{a^{2}+b^{2}}
$. \ Hence the new ratio is \ \ \ 
\[
\lim_{b\rightarrow 0^{+}}\frac{\sqrt{a^{2}+b^{2}}+\dfrac{\pi }{2}%
+\dint\limits_{0}^{1}\sqrt{\dfrac{(1-av)^{2}+b^{2}}{1-v^{2}}}\,dv}{\sqrt{%
a^{2}+b^{2}}+\dfrac{\pi }{2}+\dint\limits_{-1}^{0}\sqrt{\dfrac{%
(1-av)^{2}+b^{2}}{1-v^{2}}}\,dv}=\frac{2\sqrt{a^{2}-1}+2\func{arccsc}(a)}{%
\pi +2a}. 
\]%
After differentiating, the new equation to be solved is%
\[
\frac{\pi }{2a}\sqrt{a^{2}-1}=\func{arccsc}(a) 
\]%
which is identical to the old one in Section 4.1; therefore%
\[
a=1.2437608987462040336147443... 
\]%
reappears here. \ The objective function is different, however, and the
infimum of ratios is $0.5946...$.

\begin{tabular}{lll}
& Steven R. Finch &  \\ 
& Dept. of Statistics &  \\ 
& Harvard University &  \\ 
& Cambridge, MA, USA &  \\ 
& \textit{steven\_finch@harvard.edu} & 
\end{tabular}

\end{document}